\input amstex
\documentstyle{amsppt}
\magnification=\magstep1
\TagsOnRight
\hsize=5.1in                                                  
\vsize=7.8in
\define\R{{\bold R}}

\vskip 2cm
\topmatter

\title von Neumann Betti numbers and Novikov type inequalities\endtitle

\rightheadtext{von Neumann Betti numbers and Novikov type inequalities}
\leftheadtext{Michael Farber}
\author  Michael Farber\endauthor
\address
School of Mathematical Sciences,
Tel-Aviv University,
Ramat-Aviv 69978, Israel
\endaddress
\email farber\@math.tau.ac.il
\endemail
\thanks{\footnotemark Partially supported by US - Israel Binational Science Foundation, by
the Herman Minkowski Center for Geometry, and by EPSRC grant GR/M20563} 
\endthanks
\abstract{In this paper we show that Novikov type inequalities for closed 1-forms hold with the
von Neumann Betti numbers replacing the Novikov numbers. As a consequence we obtain a
vanishing theorem for $L^2$ cohomology. We also prove that von Neumann
Betti numbers coincide with the Novikov numbers for free abelian coverings. 
}
\endabstract
\endtopmatter

\define\C{{\bold C}}

\define\Z{{\bold Z}}

\define\V{{\Cal V}}

\define\T{{\Cal T}}


\define\im{\operatorname{im}}


\define\rank{\operatorname{rank}}

\def\<{\langle}
\def\>{\rangle}

\define\pd#1#2{\dfrac{\partial#1}{\partial#2}}

\redefine\L{\Cal L}

\def\part{\partial}

\NoBlackBoxes

\def\m{\frak m}

\def\ml{{\operatorname{Mon_L}}}
\define\n{\frak n}

\heading{\bf \S 0. Introduction}\endheading

S. Novikov and M. Shubin \cite{NS} proved that Morse inequalities for smooth functions
remain true with the usual
Betti numbers being replaced by the von Neumann Betti numbers. 

Novikov \cite{N} initiated an analog of the Morse theory for closed 1-forms. He showed that the Morse inequalities for functions can be generalized to closed 1-forms if instead of the Betti
numbers one uses {\it the Novikov numbers} (they will be briefly reviewed in \S 3.1 below).

In this paper we show that the inequalities of Novikov and Shubin \cite{NS} hold for closed 
1-forms as well. This result
gives new Novikov type inequalities for closed 1-forms. Viewed differently, this gives a vanishing theorem for $L^2$ cohomology, generalizing a theorem of W. L\"uck \cite{L3}. 

The proof of our Theorem 1 uses an idea of Gromov and Eliashberg \cite{EG}, which consist in 
counting additional critical
points which appear when transforming the given closed 1-form into a function;
the proof also uses the multiplicativity of the von Neumann Betti 
numbers under finitely sheeted coverings.

In papers \cite{BF2} and \cite{MS} different Novikov type inequalities using the von Neumann Betti
numbers were put forward. These inequalities involve 
cohomology of flat bundles of Hilbertian modules twisted
by a {\it generic} line bundle determined by the closed 1-form (similarly to the finite dimensional case).
Our approach in this paper does not require the twisting, and therefore it is simpler. On the other hand 
V. Mathai and M. Shubin \cite{MS} allow more general Hilbertian
flat bundles.

It is clear that one may easily generalize Theorem 1 below for closed 1-forms with Bott type 
singularities, in a fashion similar to \cite{BF1}, \cite{BF2}. 

I would like to thank Andrew Ranicki for stimulating discussions.

\heading{\bf \S 1. The main results}\endheading

First we recall the basic notions.

A closed 1-form $\omega$ locally is a differential of a function $\omega|_U = df_U$, where
$f_U$ is a smooth function on $U$, determined up to a constant. 
Hence, all local properties of 
functions
(e.g. the notions of a critical point, Morse singularities, indices) immediately generalize to closed
1-forms. For example, a point $p\in U$ is {\it a critical point of $\omega$ if $p$} if it is a critical point 
of $f_U$ (which is equivalent to the requirement that $\omega_p=0$); $p$ {\it a nondegenerate (or Morse) critical point } of the form $\omega$ if it is
a Morse critical point of the function $f_U$, etc.

\proclaim{1. Theorem} Let $M$ be a closed
smooth manifold and let $\omega$ be a closed 1-form on $M$,
having only Morse singularities. Let $\pi'\subset \pi=\pi_1(M)$ 
be a normal subgroup with the following property:
for any loop $\gamma \in \pi'$ holds
$\int_\gamma \omega = 0.$
Then 
$$\sum_{k=0}^i (-1)^k c_{i-k}(\omega) \, \ge \, \sum_{k=0}^i (-1)^k b_{i-k}^{(2)}(\tilde M;\pi/\pi'),\quad
i=0, 1, 2, \dots,\tag1-1$$
where $c_i(\omega)$ denotes the number of critical points of $\omega$ having index $i$, $\tilde M$
denotes the covering of $M$ corresponding to $\pi'$, and  $b^{(2)}_i(\tilde M;\pi/\pi')$ denotes the von 
Neumann Betti number of $\tilde M$. In particular, (1-1) implies
$$c_i(\omega) \, \ge\,  b^{(2)}_i(\tilde M;\pi/\pi').\tag1-2$$
\endproclaim

For the definitions of von Neumann Betti numbers and their main properties 
we refer to \cite{A}, \cite{CG}, \cite{L2}, \cite{F2}.

Theorem 1 obviously implies Corollaries 2 and 3, cf. below. These Corollaries
generalize a theorem of Wolfgang L\"uck \cite{L1} (conjectured by M. Gromov),
stating that {\it $L^2$-Betti numbers of a manifold, fibering over the circle, vanish.}

\proclaim{2. Corollary (Vanishing Theorem)} Let $M$ be a closed smooth manifold admitting a closed
1-form with Morse type singularities having no critical points of index $j$ for some $0<j<n$. Then for
any normal subgroup $\pi'\subset \pi$, lying in the commutator subgroups of $\pi=\pi_1(M)$, 
the $j$-dimensional von Neumann
Betti number vanishes
$$b^{(2)}_j(\tilde M; \pi/\pi')=0.\tag1-3$$
\endproclaim

\proclaim{3. Corollary (Morse Lacunary Principle)} Let $M$ be a closed manifold admitting a closed
1-form $\omega$ with Morse type singularities so that all the critical points of $\omega$ have even indices. Then for
any normal subgroup $\pi'\subset \pi$, lying in the commutator subgroups of $\pi=\pi_1(M)$, the 
odd-dimensional von Neumann Betti numbers vanish 
$b^{(2)}_{2j-1}(\tilde M; \pi/\pi')=0$ and the 
even-dimensional von Neumann Betti numbers are given by the formula
$$b^{(2)}_{2j}(\tilde M; \pi/\pi')=c_{2j}(\omega).\tag1-4$$
\endproclaim
Note that Theorem 0.3 of \cite{BF1} implies that under the conditions of Corollary 2 the number
$c_{2j}(\omega)$ equals the Novikov number $b_{2j}(\xi)$, where $\xi\in H^1(M;\R)$ is the cohomology
class of $\omega$. Hence we obtain the equality
$$b^{(2)}_{2j}(\tilde M; \pi/\pi') = b_{2j}(\xi),$$ 
compare Theorem 5 below.

Definitions of the Novikov numbers can be found in the literature
\cite{BF1}, \cite{F1}, \cite{N}. Paper \cite{BF1} contains a stronger version of Novikov type inequalities and a more general notion of Novikov numbers. 
In this paper we will use only the standard Novikov numbers $b_i(\xi)$. It will be convenient for our
purposes to be based on the definition given in \cite{BF1}, cf. \S 3.1.

\proclaim{4. Corollary (Symplectic circle actions)} 
Let $M$ be a closed symplectic manifold admitting a symplectic circle
action with isolated fixed points. Then for
any normal subgroup $\pi'\subset \pi$, lying in the commutator subgroups of $\pi=\pi_1(M)$, the 
odd-dimensional von Neumann Betti numbers vanish 
$b^{(2)}_{2j-1}(\tilde M; \pi/\pi')=0$, and the 
even-dimensional von Neumann Betti numbers coincide with the Novikov numbers
$$b^{(2)}_{2j}(\tilde M; \pi/\pi')=b_{2j}(\xi).\tag1-5$$
Here $\xi\in H^1(M;\R)$ denotes the cohomology class of the generalized moment map, and 
$b_i(\xi)$ denotes the $i$-dimensional Novikov number, corresponding to $\xi$. 
\endproclaim

\demo{Proof of Corollary 4} Let us first explain the terms used in Corollary 4.
Suppose that $\Omega$ denotes the symplectic form
of $M$. The $S^1$ action is assumed to be {\it symplectic}, which means that for any $g\in S^1$ holds 
$g^\ast\Omega =\Omega$. Let $X$ denote the vector field
generating the $S^1$-action. Then 
$$\omega = \iota(X)\Omega\tag1-6$$
is a closed 1-form on $M$, which
is called the {\it generalized moment map}. We consider its De Rham cohomology class 
$\xi=[\omega] \in H^1(M;\R)$ of $\omega$. 
Recall that a symplectic circle action is called {\it Hamiltonian} if $\xi=0$. In this case the Novikov numbers $b_i(\xi)$ coincide
with the Betti numbers $b_i(M)$.

The critical points of the generalized moment map $\omega$ are precisely the fixed points of the 
circle action. It is well known that $\omega$ has Morse type singularities,
assuming that the fixed points are isolated; moreover,
all the critical points have even indices, cf. \cite{Au}. Hence from Corollary 3 we obtain
that $b^{(2)}_{2j-1}(\tilde M; \pi/\pi')=0$ and the 
even-dimensional von Neumann Betti numbers are given by the formula (1-4).
Now we will use the main theorem 0.3 of \cite{BF1}, which gives
$$c_{2j}(\omega)=b_{2j}(\xi).\tag1-7$$
\qed
\enddemo

The following Theorem states that in the special case when $\pi'$ coincides with the $\ker(\xi)$
the inqualities of Theorem 1 are precisely the Novikov inequalities.

\proclaim{5. Theorem (Novikov numbers equal $L^2$ Betti numbers)} Let $M$ be a closed
smooth manifold and let $\xi\in H^1(M;\R)$ be a nontrivial cohomology class. 
Let $\pi'$ be the kernel of the homomorphism $\xi:\pi=\pi_1(M) \to \R$ determined
by $\xi$. Then the $i$-dimensional Novikov  number $b_i(\xi)$ coincides with the $i$-dimensional
von Neumann Betti number of the covering $\tilde M \to M$ corresponding to $\pi'$:
$$b_i(\xi) \, =\,  b_i^{(2)}(\tilde M; \pi/\pi'). \tag1-8$$
\endproclaim

The plan of the paper is as follows.
Proof of Theorem 1 is presented in \S 2. Proof of Theorem 5 is described in \S 3.

\heading{\bf \S 2. Proof of Theorem 1}\endheading

Let $\omega$ be a closed 1-form on a closed manifold $M$.
First, we will assume that the form $\omega$ has integral periods,
i.e. its cohomology class is integral $\xi=[\omega]\in H^1(M;\Z)$; the general case will be treated at the
end of the proof.

There exists a smooth map into the circle $f: M\to S^1$ with $\omega = f^\ast(d\theta)$,
where $d\theta$ denotes the standard angular form on the circle. This map is given by
$$f(x) = \exp(2\pi i\int_{x_0}^x \omega),$$
where $x_0\in M$ is a base point. We will also denote by $\xi:\pi=\pi_1(M) \to \Z$
the induced homomorphism $[\gamma]\mapsto \int_\gamma \omega$, where $\gamma$ is a loop in
$M$.

Given a positive integer $m$, we will denote by $\pi_m$ the preimage $\pi_m =\xi^{-1}(m\Z)$. 
It is a normal subgroup containing $\pi'$. Denote by $\tilde M_m$ the cyclic $m$-sheeted
covering $\tilde M_m\to M$ corresponding to the subgroup $\pi_m$. 
It is clear that the pullback $\tilde \omega_m$ of the form $\omega$ to $\tilde M_m$ has
the following Morse numbers
$$c_j(\tilde \omega_m) = m\cdot c_j(\omega).\tag2-1$$

\subheading{Claim} {\it There exists a constant $C$,  independent of $m$, so that on each 
manifold $\tilde M_m$ there exists an exact 1-form $\tilde \omega'_m$}
$$c_j(\tilde \omega_m) \le c_j(\tilde \omega'_m) \le c_j(\tilde \omega_m) + C, 
\qquad j=0,1, \dots\tag2-2$$

We will postpone the proof of the Claim and will continue the proof of the Theorem.

We have the inequality
$$\sum_{k=0}^i (-1)^k c_{i-k}(\tilde \omega'_m) \, \ge \, 
\sum_{k=0}^i (-1)^k b_{i-k}^{(2)}(\tilde M;\pi_m/\pi'),\quad
i=0, 1, 2, \dots,\tag2-3$$
which is just 
the Novikov - Shubin inequality \cite{NS} (cf. also \cite{MS}) applied to the exact form $\tilde \omega'_m$ on the compact manifold $\tilde M_m$.

Multiplicativity of the von Neumann Betti numbers under finitely sheeted coverings gives
$$b_{j}^{(2)}(\tilde M;\pi_m/\pi') \, = \, m \cdot b_{j}^{(2)}(\tilde M;\pi/\pi'),\tag2-4$$
cf. \cite{L2}, Theorem 1.7, statement 7.
Here we use our assumption that $\pi'\subset \ker(\xi)$.

Combining (2-1), (2-2), (2-3), (2-4)  we will have
$$
\aligned
&\sum_{k=0}^i (-1)^k c_{i-k}(\omega) \, = \, \\
& \frac{1}{m} \sum_{k=0}^i (-1)^k c_{i-k}(\tilde \omega_m) \, \ge \,\\
& \frac{1}{m} [\sum_{k=0}^i (-1)^k c_{i-k}(\tilde \omega'_m) - (i+1)\cdot C]\, \ge \,\\
& \frac{1}{m}[ \sum_{k=0}^i (-1)^k b^{(2)}_{i-k}(\tilde M;\pi_m/\pi') - (i+1)\cdot C] =\\
&  \sum_{k=0}^i (-1)^k b^{(2)}_{i-k}(\tilde M;\pi/\pi') - \frac{(i+1)\cdot C}{m}.
\endaligned\tag2-5
$$
Hence,
$$ \sum_{k=0}^i (-1)^k c_{i-k}(\omega) \, \ge
 \sum_{k=0}^i (-1)^k b^{(2)}_{i-k}(\tilde M;\pi/\pi') - \frac{(i+1)\cdot C}{m}\tag2-6
$$
for arbitrary $m>0$. 
Since $C$ is independent of $m$, taking the limit $m\to \infty$ in the obtained inequality (2-6)
proves (1-1). 

Now we want to prove the Claim.

Suppose that $\exp(2\pi i \theta_0)\in S^1$ is a regular 
value of the map $f: M\to S^1$, where $0<\theta_0<1$. 
Let $V\subset M$ denote $f^{-1}(\theta_0)$.
Consider the cylinder $V\times [-1,1]$ and find an arbitrary Morse function 
$$\phi: V\times [-1,1]\to \R\tag2-7$$ 
with the following properties:

(a) $\phi$ assumes values in the open interval $(0,1)$;

(b) $\phi(v,t) =1+\delta t$ for all $t\in [-1,-1/2)$ and $v\in V$, where $0<\delta <1$ is a fixed number;

(c) $\phi(v,t) =\delta t$ for all $t\in (1/2, 1]$ and $v\in V$.

We will denote by $C$ the total number of critical points of $\phi$. 

For any $m$ we have the following commutative diagram
$$
\CD
\tilde M_m @>{f_m}>> S^1\\
@VVV                          @VV{z\mapsto z^m}V\\
M_m @>f>> S^1
\endCD
$$
We obtain
$\tilde \omega_m \, = f_m^\ast(d\theta),$
where $f_m : \tilde M_m \to S^1$ is a smooth map.
It is clear that for any $m$ we may smoothly imbed $V\times [-1,1]$ into $\tilde M_m$ so that
on the image of $V\times [-1,1]$ the map $f_m$ is given by
$f_m(v,t) = \exp(2\pi i[\theta_0/m +\delta t]).$
We will construct now a new map $g_m : \tilde M_m \to S^1$ as follows:

(a) $g_m$ coincides with $f_m$ on the complement of $V\times [-1,1]$ in $\tilde M_m$;

(b) on $V\times [-1,1]$ the map $g_m$ is given by
$$g_m(v,t) = \exp(2\pi i[\theta_0/m +\phi(v,t)]),$$
where $\phi: V\times [-1,1]\to \R$ is the function constructed above.

Set $\tilde \omega'_m = g_m^\ast(d\theta)$. Then the obtained 1-form $\tilde \omega'_m$ is 
{\it exact} since the
map $g_m: \tilde M_m\to S^1$ is null homotopic (it does not assume the 
value $\exp(2\pi i \theta_0/m)\in S^1$). The inequality (2-2) is clearly satisfied since 
the forms $\tilde\omega_m$ and $\tilde\omega'_m$ have the same critical points on the complement 
of $V\times [-1,1]$ and on $V\times [-1,1]$ the form $\tilde\omega_m$ has no critical points, while
the form $\tilde\omega'_m$ may have at most $C$ critical points. 

This proves the Claim and hence completes the 
proof of Theorem 1 in the case of closed 1-forms $\omega$ having integral cohomology classes.

Consider now a closed 1-form $\omega$ representing an arbitrary cohomology class 
$\xi=[\omega]\in H^1(M;\R)$ and having Morse type singularities. It clearly defines a homomorphism
$\xi: H_1(M;\Z)\to \R$ and the image $\im(\xi)$ 
is a finitely generated free abelian group. We will denote the
rank of this image by $r=\rank(\im(\xi)) = \rank(\omega)$. 

Note that any closed 1-form $\omega$ having rank 1 
can be represented as $\omega=\lambda\omega_0$, where $\lambda>0$ and $\omega_0$ represents
an integral class. Since the forms $\omega$ and $\omega_0$ have the same critical points, we 
conclude
that {\it Theorem 1 holds for all closed 1-forms of rank 1.}

Suppose now that $\omega$ is a closed 1-form having rank $r>1$. Then we may find a sequence
of closed 1-forms $\omega_n$ on $M$ with the following properties:

(i) $\omega_n$ converges to $\omega$ in the $C^0$-norm;

(ii) there exists a neighbourhood $U$ of the set $S$ of critical points of $\omega$ such that 
$\omega|_U = \omega_n|_U$ for all $n$;

(iii) $\omega_n$ has rank 1 for every $n$;

(iv) the homomorphism $[\omega_n]: H_1(M;\Z)\to \R$ vanishes on $\ker(\xi)$.

Because of (i) and (ii),
for large $n$ the form $\omega_n$ has the same critical points as $\omega$ and they have the
same indices $c_j(\omega) = c_j(\omega_n)$.
Because of properties (iii) and (iv) we have the inequality 
$$\sum_{k=0}^i (-1)^k c_{i-k}(\omega_n) \, \ge \, 
\sum_{k=0}^i (-1)^k b_{i-k}^{(2)}(\tilde M;\pi/\pi'),\quad
i=0, 1, 2, \dots,$$
which obviously implies (1-1). 

Let us show how to construct the sequence of forms $\omega_n$. Let $N_\xi\subset H^1(M;\R)$
be the subspace formed by classes $\eta$ with $\eta|_{\ker(\xi)}=0$; here $r=\rank(\xi)$. Let
$\eta_1, \dots, \eta_r$ be a basis of $N_\xi$ such that $\rank(\eta_j)=1$ for $j=1, \dots, r$. We may
realize each class $\eta_j$ by a closed 1-form $\nu_j$ which vanishes identically in a neighbourhood
of the zeros of $\omega$. We may write $\xi=\sum_{j=1}^r\alpha_j\eta_j$, where $\alpha_j\in \R$.
Now, choose a sequence of rational numbers $\alpha^n_j$ converging to $\alpha_j$ as $n\to \infty$.
Then the sequence of closed 1-forms
$$\omega_n = \omega + \sum_{j=1}^r (\alpha_j^n-\alpha_j)\nu_j$$
satisfies all the requirements.
This completes the proof. \qed

\heading{\bf \S 3. Proof of Theorem 5}\endheading

The proof of theorem 5 will use basic harmonic analysis. 

Before starting the proof we recall the definition of the
Novikov numbers.

\subheading{3.1. The Novikov numbers $b_i(\xi)$} Let $\xi\in H^1(M;\R)$ be a cohomology
class. We will denote by $\ker(\xi)$ the kernel of the induced homomorphism $\xi: H_1(M;\Z)\to \R$.

Let us consider complex flat line bundles $L \to M$ with the following property: the 
monodromy along
any loop $\gamma\in \ker(\xi)$ is identity. We will denote by $\V_\xi$ the set of isomorphism classes
of all such flat line bundles. Given $L\in \V_\xi$, it has the monodromy homomorphism
$$\ml: H_1(M;\Z)/\ker(\xi) \to \C^\ast,$$
which completely determines the isomorphism type of $L$. Hence, there is a one-to-one correspondence between $\V_\xi$ and the complex torus $(\C^\ast)^r$. Here $r$ is the rank of class $\xi$, i.e. $r=\rank (H_1(M;\Z)/\ker(\xi)) =\rank(\xi)$. This identification $\V_\xi \simeq (\C^\ast)^r$
allows to view $\V_\xi$ as an affine algebraic variety.

Fix some $i$ (with $0\le i\le \dim M)$ and consider the function
$$\V_\xi \ni L \, \, \mapsto\, \, \dim_\C H_i(M;L).\tag3-1$$
It is well know in the algebraic geometry (cf. \cite{H}, chapter 3, \S 12) that function (3-1)
has the following property: {\it there exists a proper algebraic subvariety 
$V=V_i(M)\subset \V_\xi$ such that
the dimension $\dim_\C H_i(M;L)$ is constant for all $L\notin V$ and for $L'\in V$}
$$\dim_\C H_i(M;L') > \dim_\C H_i(M;L).\tag3-2$$ 

\subheading{Definition} {\it The Novikov number $b_i(\xi)$ is defined as 
$\dim_\C H_i(M;L)$ for $L\in \V_\xi - V$.}

\subheading{3.2. Unitary flat bundles} Consider the subset $\T_\xi\subset \V_\xi$ 
consisting of isomorphism classes of flat
line bundles $L$ admitting flat Hermitian metrics. Under the identification $\V_\xi \simeq (\C^\ast)^r$
given by the monodromy representation, the subset $\T_\xi$ corresponds to the real torus
$(S^1)^r\subset (\C^\ast)^r$. This shows, in particular, 
that $\T_\xi$ is a real analytic subvariety of 
$\V_\xi$.

It is easy to see that {\it 
any complex Laurent polynomial $p(z_1, \dots, z_r, z_1^{-1}, \dots, z_r^{-1})$, which vanishes
on the torus $(S^1)^r\subset (\C^\ast)^r$, is identically zero.} Hence, the intersection of any non-empty Zariski
open subset of $(\C^\ast)^r$ with the torus $\T_\xi$ is non-empty. 

This implies the following interpretation
of the Novikov numbers:

\proclaim{3.3. Corollary} For a fixed $M$ and $\xi\in H^1(M;\R)$, consider the function 
$$T_\xi \ni L \, \, \mapsto\, \, \dim_\C H_i(M;L).\tag3-3$$
There exists a proper real analytic subvariety $V=V_i(M)\subset \T_\xi$ such that for all 
$L\notin V$ 
the dimension $\dim_\C H_i(M;L)$ equals to the Novikov number $b_i(\xi)$
and for $L'\in V$
$$\dim_\C H_i(M;L') > b_i(\xi).\tag3-4$$ 
\endproclaim

\subheading{3.4. Proof of Theorem 5} Let $M$ be a closed manifold and let $\xi\in H^1(M;\R)$ be a real cohomology class. Consider the covering $\tilde M \to M$ corresponding to the subgroup
$\pi'=\ker(\xi)\subset \pi=\pi_1(M)$. Then $\pi/\pi'$ is isomorphic to $\Z^r$ for some $r$.

Recall that the von Neumann Betti number 
$b_i^{(2)}(\tilde M;\pi/\pi')$ is defined as follows. 
Fix a smooth triangulation of $M$ and consider the
induced (equivariant) triangulation of the covering space $\tilde M$. Let $C_\ast(\tilde M)$ denote the
simplicial chain complex of this triangulation; it consists of free finitely generated $\Z[\Z^r]$-modules
and $\Z[\Z^r]$-module homomorphisms. Consider the Hilbert space $L^2(\Z^r)$ with its canonical
$\Z^r$-action and form the complex of Hilbert spaces 
$$L^2(\Z^r)\otimes_{\Z^r} C_\ast(\tilde M).\tag3-5$$
The von Neumann Betti number 
$b_i^{(2)}(\tilde M;\pi/\pi')$ is defined as the von Neumann dimension of the reduced homology of
complex (3-5), cf. \cite{A}, \cite{CG}, \cite{F2}, \cite{L2}.

Let us now remind the standard construction of the harmonic analysis, providing
an isomorphism between the space of $L^2$-functions on the torus $(S^1)^r$ with respect to the
Lesbegue measure $d\mu$ (which will be assumed to be normalized so that the total torus 
$(S^1)^r$ has measure 1) and the space
of $L^2$-functions on the lattice $\Z^r$. Points of the lattice $\Z^r$ we will denote by 
$r$-tuples
$\n=(n_1, n_2, \dots, n_r)$, where $n_1, \dots, n_r\in \Z$. 
Points of the torus $(S^1)^r$ we will denote by $z=(z_1, z_2, \dots, z_r)$, 
where $z_1, \dots, z_r\in S^1\subset \C$. The symbol $z^\n$ will denote the product
$z^{\n}= z_1^{n_1}\cdot z_2^{n_2}\cdot \dots \cdot z_r^{n_r}\in S^1.$
Any $L^2$-function $f\in L^2((S^1)^r,d\mu)$ can be uniquely represented by its Fourier series
$$f(z) = \sum_{\n\in \Z^r} a_\n z^\n.$$
The correspondence $f\mapsto (a_\n)_{\n \in \Z^r}$ defines an isometry 
$L^2((S^1)^r,d\mu)\simeq L^2(\Z^r)$. It is important that under this isometry the multiplication
by $z^\m$ in $L^2((S^1)^r,d\mu)$ (i.e. the map $f(z)\mapsto z^\m f(z)$) transforms into the shift
$(a_\n)_{\n\in \Z^r}\mapsto (a_\n)_{\n+\m\in \Z^r}$ in $L^2(\Z^r)$; here $\m\in \Z^r$.

The above identification allows to view the Hilbert space $L^2(\Z^r)$ as the space of $L^2$-sections
of a vector bundle over the torus $\T_\xi$, which we will now describe.
Consider the real analytic line bundle $\L\to \T_\xi$ of unitary $\Z^r$-module structures
on $\C$. It is trivial as a vector bundle and the fiber over a point $z\in \T_\xi$ has the following
$\Z^r$-action: a point of the lattice $\n\in \Z^r$ acts as multiplication by $z^\n\in S^1$. We may form
the chain complex
$$\L\otimes_{\Z^r} C_\ast(\tilde M)\tag3-6$$
which is a real analytic family (parametrized by the points of the torus $\T_\xi$) of chain 
complexes of finite dimensional
vector spaces. The $L^2$-complex (3-5) is isomorphic to the complex of $L^2$-sections of (3-6).
Hence, applying Theorem 4.11 from \cite{F3} we obtain that the von Neumann Betti number
$b_i^{(2)}(\tilde M; \pi/\pi')$ coincides with the {\it generic Betti number}
$\dim_\C H_i(M;L)$ for $L\in \T_\xi -V$, where $V$ is a proper real analytic subvariety.
This may be also stated as the equality
$$b_i^{(2)}(\tilde M; \pi/\pi') = \int_{\T_\xi} \dim_\C H_i(M;L)d\mu(L).\tag3-7$$
By Corollary 3.3 the generic Betti number is precisely the Novikov number $b_i(\xi)$. 
This completes the proof. \qed

\enddocument